\documentclass[11pt]{amsart}
\usepackage[margin=1.1in]{geometry}
\usepackage{amsmath,amssymb,amsthm,mathtools}
\usepackage{graphicx}
\usepackage{booktabs}
\usepackage[labelfont=bf,labelsep=space]{caption}
\usepackage{placeins}
\usepackage{url}
\usepackage{cite}
\renewcommand{\texttt}[1]{\nolinkurl{#1}}
\usepackage[colorlinks=true,linkcolor=blue,citecolor=blue,urlcolor=blue]{hyperref}

\newtheorem{theorem}{Theorem}[section]
\newtheorem{proposition}[theorem]{Proposition}
\newtheorem{lemma}[theorem]{Lemma}
\newtheorem{corollary}[theorem]{Corollary}
\theoremstyle{definition}

\newtheorem{assumption}[theorem]{Assumption}
\theoremstyle{remark}
\newtheorem{remark}[theorem]{Remark}

\numberwithin{equation}{section}
\newcommand{\N}{\mathbb{N}}
\newcommand{\R}{\mathbb{R}}
\newcommand{\E}{\mathbb{E}}
\newcommand{\PP}{\mathbb{P}}

\newcommand{\eps}{\varepsilon}
\newcommand{\rhoq}{\rho_{\mathrm{q}}}
\newcommand{\rhoa}{\rho_{\mathrm{a}}}
\newcommand{\cF}{\mathcal{F}}
\newcommand{\cP}{\mathcal{P}}
\newcommand{\cS}{\mathcal{S}}
\newcommand{\Imin}{I_{-}}
\newcommand{\Ip}{I_{+}}
\newcommand{\Prw}{\mathbf{P}}
\newcommand{\Erw}{\mathbf{E}}
\DeclareMathOperator*{\esssup}{ess\,sup}
\DeclareMathOperator*{\essinf}{ess\,inf}

\title[Lower deviations for the maximum of a BRW in random environment in time]{Lower deviations for the maximum of a branching random walk in a random environment in time: quenched versus annealed}

\author[W. Hou]{Wanting Hou}
\thanks{Corresponding author: Wanting Hou, School of Mathematics and Statistics, Liaoning University, Shenyang 110801, P. R. China. E-mail: \href{mailto:houwanting@lnu.edu.cn}{\nolinkurl{houwanting@lnu.edu.cn}}.}

\begin{document}

\begin{abstract}
We study lower deviations of the rightmost position in a supercritical branching random walk whose offspring laws are independent and identically distributed across generations. Every particle has at least one child, and displacements have a fixed law with finite exponential moments of all orders. For thresholds strictly between the lower endpoint of the step law and the speed of the maximum, the quenched lower deviation probability has a finite positive exponential rate precisely when the negative logarithm of the probability of a single child is integrable. We identify the rate by a variational formula. If a positive fraction of generations rule out single children and displacements are bounded below, the probability decays doubly exponentially. We compute its double-exponential rate in terms of the minimum offspring numbers, the lower endpoint of the displacement law, and the speed of the maximum. After averaging over the environment, the negative logarithm of the probability is bounded above and below by positive multiples of the generation number whenever single children are possible with positive probability. We obtain the exact exponential rate when the offspring mean is constant and the offspring second moments are uniformly bounded. With displacements bounded below, mixed environments thus yield double-exponential quenched decay and single-exponential annealed decay.
\end{abstract}


\maketitle

\section{Introduction}\label{sec:intro}

\subsection{Background}
Let $M_n$ denote the maximal position at generation $n$ of a supercritical branching random walk (BRW) on $\R$. By the results of Hammersley \cite{Hammersley1974}, Kingman \cite{Kingman1975} and Biggins \cite{Biggins1977}, $M_n/n$ converges almost surely on survival to a constant $x^*$. The logarithmic correction was identified by Hu and Shi \cite{HuShi2009} and by Addario-Berry and Reed \cite{AddarioBerryReed2009}, and A\"id\'ekon \cite{Aidekon2013} proved, under a non-lattice assumption, that $M_n-x^*n+c\log n$ converges in law for an explicit $c>0$; see \cite{Shi2015} for an account. The probability that $M_n$ is atypically small,
\[
\PP(M_n\le xn),\qquad x<x^*,
\]
behaves very differently from the upper deviation probability $\PP(M_n\ge xn)$, $x>x^*$. Consider a constant environment, i.e.\ a Galton--Watson tree with offspring law $(p_k)_{k\ge0}$ and first generation size $Z_1\sim(p_k)$. Gantert and H\"ofelsauer \cite[Theorem~3.2]{GantertHofelsauer2018} showed that in the \emph{Schr\"oder case} $p_0+p_1>0$ (the terminology is standard in the lower deviation theory of Galton--Watson processes, see e.g.\ \cite{NeyVidyashankar2003,FleischmannWachtel2007}) the lower deviation probability (conditioned on survival if $p_0>0$) decays exponentially, with rate
\begin{equation}\label{eq:GH}
H(x)=\inf_{t\in(0,1]}\Bigl\{t\rho+t\,I\Bigl(\frac{x-(1-t)x^*}{t}\Bigr)\Bigr\},
\end{equation}
where $I$ is the Cram\'er rate function of the step distribution, $\rho=-\log\E[Z_1q^{Z_1-1}]$ (expectation under the Galton--Watson law) and $q$ is the extinction probability; in particular $\rho=-\log p_1$ when $p_0=0$. The optimal strategy keeps a single particle alive during the first $tn$ generations and lets it move atypically low. Chen and He \cite[Theorem~1.1]{ChenHe2020} treated the \emph{B\"ottcher case} $p_0=p_1=0$, where every particle has at least $b:=\min\{k:p_k>0\}\ge2$ children. There the lower deviation probability decays doubly exponentially; if the steps are bounded from below by $-L$, then
\begin{equation}\label{eq:CH}
\lim_{n\to\infty}\frac1n\log\bigl[-\log\PP(M_n\le xn)\bigr]=\frac{x^*-x}{x^*+L}\log b,\qquad x\in(-L,x^*).
\end{equation}
Derrida and Shi \cite{DerridaShi2017} obtained analogous results for branching Brownian motion and some branching random walks (see also \cite{Hu2016}), and the regime between typical and large deviations was studied in \cite{ChenHe2020} (lower side) and in the preprint by Chataignier and Luo (\href{https://arxiv.org/abs/2601.08766}{arXiv:2601.08766}; upper side). Thus, in a constant environment and for steps with finite exponential moments of all orders, whether the lower deviation probability (conditioned on survival if $p_0>0$) decays exponentially or faster is determined by whether $p_0+p_1>0$ (for heavier left tails see \cite[Theorem~1.3]{ChenHe2020}). Lower deviation probabilities for the level sets of the branching random walk, that is, for the numbers of particles above a linear level, were studied by Zhang \cite{Zhang2023}.

A branching random walk in a random environment in time (BRWre) is a BRW in which the reproduction and displacement laws of generation $k$ are given by the $k$-th coordinate $\xi_k$ of an i.i.d.\ sequence. Its typical behaviour is by now well understood: the first and second order of $M_n$ were identified by Mallein and Mi{\l}o{\'s} \cite{MalleinMilos2019} (see also \cite{FangZeitouni2012,Mallein2015,HouZhangHong2021} for deterministic time-inhomogeneous environments), tightness was studied by Kriechbaum \cite{Kriechbaum2025}, the derivative martingale by Hong and Liang \cite{HongLiang2025}, limit theorems for the counting measures by Huang and Liu \cite{HuangLiu2024}, the case of critical branching ($\E\log m=0$) by Fu and Hong \cite{FuHong2026}, heavy-tailed displacements by Bhattacharya and Palmowski \cite{BhattacharyaPalmowski2025}, and the extinction rate of a BRWre with a barrier by Lv \cite{Lv2025}; large and moderate deviations of the counting measures were obtained in \cite{HuangLiangLiu2014,WangHuang2017}.

Rare events in which the population stays small have been analysed for branching processes in random environment: lower large deviations and small positive values of $Z_n$ \cite{BansayeBerestycki2009,HuangLiu2012,BansayeBoinghoff2013,BansayeBoinghoff2014} (see \cite{KerstingVatutin2017} for background). In these works the quenched and annealed costs of keeping a population with $p_0=0$ bounded are $\E[-\log p_1]$ and $-\log\E[p_1]$; what is new here is the interplay of this cost with the displacement of the surviving particles and with the speed in the shifted environment. Quenched small deviations of a single random-walk path in a random environment in time were studied in \cite{LvHong2023}. To our knowledge, the lower deviations of the maximal displacement $M_n$ of a BRWre have not been studied. In this paper the step law does not depend on the environment (unlike in \cite{MalleinMilos2019}); only the offspring laws are random.

\subsection{The question and the answer in brief}
A random environment can mix generations of Schr\"oder type ($p_1(\xi_k)>0$) and of B\"ottcher type ($p_1(\xi_k)=0$). This paper answers the following question: \emph{how does the environment determine the decay of $\PP(M_n\le xn)$, $x<x^*$, under the quenched and under the annealed law?} The answer has three parts.

\begin{enumerate}
\item \emph{Quenched, $\rhoq<\infty$} (Theorem~\ref{thm:A}). If $\rhoq=\E[-\log p_1(\xi_1)]<\infty$, then for almost every environment $P_\xi(M_n\le xn)=e^{-nH(x;\rhoq)+o(n)}$, where $H(x;\rho)$ is defined in \eqref{eq:H} (it coincides with \eqref{eq:GH} for $\rho=-\log p_1$). If $\rhoq=\infty$, the decay is faster than exponential; apart from the case in item 2, we do not determine its order.
\item \emph{Quenched, B\"ottcher-type generations with positive probability} (Theorem~\ref{thm:B}). If the steps are bounded below by $-L>-\infty$, even a small proportion of B\"ottcher-type generations makes the decay doubly exponential: for almost every environment,
$-\log P_\xi(M_n\le xn)=\exp\{n\beta(x^*-x)/(x^*+L)+o(n)\}$ with $\beta=\E[\log b(\xi_1)]$.
\item \emph{Annealed} (Theorem~\ref{thm:annealed}). Whenever Schr\"oder-type generations occur with positive probability, $\PP(M_n\le xn)\ge e^{-nH(x;\rhoa)+o(n)}$ with $\rhoa=-\log\E[p_1(\xi_1)]$; if the mean number of children does not depend on the environment and the second moments of the offspring laws are uniformly bounded, this is the exact exponential rate. In general, the annealed probability decays at least exponentially fast (Proposition~\ref{prop:annealedupper}).
\end{enumerate}

Comparing these statements (Corollary~\ref{cor:dichotomy}) shows two effects that have no counterpart in a constant environment. If Schr\"oder- and B\"ottcher-type generations are mixed (and $L<\infty$), the quenched probability decays doubly exponentially while the annealed probability decays exponentially: the annealed lower bound is obtained from environments whose first $tn$ generations are all of Schr\"oder type, an event of exponentially small (and not doubly exponentially small) probability. If $\rhoq<\infty$ and $p_1(\xi_1)$ is random, the annealed probability is exponentially larger than the quenched one, because $\rhoa<\rhoq$ by Jensen's inequality.

\subsection{Method}
The lower bounds follow the strategies of \cite{GantertHofelsauer2018,ChenHe2020}, evaluated along the environment with the strong law of large numbers. The upper bounds require two estimates that hold \emph{uniformly over time-shifts of the environment}, because after the first $tn$ generations the remaining process evolves in the shifted environment $\theta^{tn}\xi$. The first (Lemma~\ref{lem:chain}) shows that the probability that the maximum reaches $x^*-\eta$ times the elapsed time is not exponentially small, uniformly over time-shifts; it is obtained by chaining blocks of length $B$ and controlling the block probabilities through an embedded branching process in random environment and the theorem of Athreya and Karlin \cite{AthreyaKarlin1971}. The second (Lemma~\ref{lem:apriori}) is an a priori exponential bound for lower deviations, uniform over shifts; with it, the upper bound in the Schr\"oder case can be organised around the first time the population reaches a fixed size $K$, which avoids estimates on population sizes growing with $n$. A crude bound on the harmonic moments of the population (Lemma~\ref{lem:harmonic}; sharper harmonic moment estimates are due to Ney and Vidyashankar \cite{NeyVidyashankar2003} for Galton--Watson processes and to Grama, Liu and Miqueu \cite{GramaLiuMiqueu2017} in random environment) provides the quenched cost of keeping the population small.

\subsection{Organisation}
Section~\ref{sec:model} introduces the model and states the results. Section~\ref{sec:prelim} collects the preliminary estimates and proves Proposition~\ref{prop:LLN}. Theorems~\ref{thm:A}, \ref{thm:B} and \ref{thm:annealed} are proved in Sections~\ref{sec:proofA}, \ref{sec:proofB} and \ref{sec:proofannealed}; Proposition~\ref{prop:annealedupper} is proved at the end of Section~\ref{sec:proofannealed}. Section~\ref{sec:numerics} reports numerical computations (deterministic for the quenched law, by enumeration or importance sampling for the annealed law), and Section~\ref{sec:discussion} discusses the assumptions.

\section{Model and main results}\label{sec:model}

\subsection{The model}
Let $\cP$ denote the set of probability measures on $\N_0=\{0,1,2,\dots\}$ (offspring laws). For $\zeta\in\cP$ write
\[
m(\zeta)=\sum_{k\ge1}k\,\zeta(k),\qquad p_1(\zeta)=\zeta(1),\qquad b(\zeta)=\min\{k\ge0:\zeta(k)>0\}.
\]
The \emph{environment} $\xi=(\xi_1,\xi_2,\dots)$ is an i.i.d.\ sequence of $\cP$-valued random variables with common law $\nu$, defined on a probability space with probability $\PP_{\mathrm{env}}$ and expectation $\E$. Let $\mu$ be a probability measure on $\R$ (the step law, the same in all generations). On a separate probability space with probability $\Prw$ and expectation $\Erw$, let $X,X_1,X_2,\dots$ be i.i.d.\ with law $\mu$ and $S_n=X_1+\dots+X_n$; this auxiliary random walk is used only through its law.

Given $\xi$, the BRW is defined as follows. At time $0$ a single particle sits at the origin. Each particle of generation $k-1$ independently produces a random number of children with law $\xi_k$, and each child is displaced from its parent by an independent step with law $\mu$. Thus the genealogical tree is a Galton--Watson tree in the varying environment $\xi$, and, given the tree, the displacements are i.i.d.\ with law $\mu$. We write $P_\xi$ and $E_\xi$ for the \emph{quenched} law and expectation, and
\[
\PP(\cdot)=\E\bigl[P_\xi(\cdot)\bigr]
\]
for the \emph{annealed} law, with expectation $\E_\PP$. Let $D_n$ be the set of particles of generation $n$, $Z_n=\#D_n$, $S_u$ the position of $u$, and
\[
M_n=\max_{u\in D_n}S_u .
\]
For $u\in D_k$ and $n\ge k$, $M^{u}_{n-k}=\max\{S_w-S_u: w\in D_n,\ w\text{ a descendant of }u\}$. Let $\cF_k$ be the $\sigma$-field generated by the first $k$ generations (genealogy and positions) and $\theta$ the shift, $(\theta\xi)_i=\xi_{i+1}$. By the branching property, conditionally on $\cF_k$ the processes $(M^u_{n-k})_{n\ge k}$, $u\in D_k$, are independent, and each has the law of $(M_{n-k})_{n\ge k}$ under $P_{\theta^k\xi}$.

\subsection{Assumptions and notation}
\begin{assumption}\label{ass:standing}
\begin{enumerate}
\item[(A1)] $\xi_1(0)=0$ almost surely, i.e.\ every particle has at least one child.
\item[(A2)] $\ell:=\E[\log m(\xi_1)]\in(0,\infty)$.
\item[(A3)] $\mu$ is not a Dirac mass and $\psi(\lambda):=\int e^{\lambda x}\mu(dx)<\infty$ for all $\lambda\in\R$.
\end{enumerate}
\end{assumption}

Under (A1)--(A2), $m(\xi_1)\ge1$ and $\PP_{\mathrm{env}}(p_1(\xi_1)<1)>0$. Let
\[
I(y)=\sup_{\lambda\in\R}\{\lambda y-\log\psi(\lambda)\},\quad
\Imin(y)=\sup_{\lambda\le0}\{\lambda y-\log\psi(\lambda)\},\quad
\Ip(y)=\sup_{\lambda\ge0}\{\lambda y-\log\psi(\lambda)\}.
\]
$\Imin$ is convex, nonincreasing, lower semicontinuous, vanishes on $[\Erw X,\infty)$ and coincides with $I$ on $(-\infty,\Erw X]$; symmetric statements hold for $\Ip$. By (A3), $\Imin(y)\to\infty$ faster than linearly as $y\to-\infty$ (if the support of $\mu$ is unbounded below) and $\Ip(y)\to\infty$ as $y\to\infty$. Let $-L:=\essinf X\in[-\infty,\infty)$, so that $L\in(-\infty,\infty]$, and define the speed
\begin{equation}\label{eq:xstar}
x^*:=\sup\{x\in\R:\Ip(x)\le\ell\}=\inf_{\lambda>0}\frac{\ell+\log\psi(\lambda)}{\lambda}.
\end{equation}
Since $\Ip(\Erw X)=0<\ell$ and $\mu$ is non-degenerate, $\Erw X<x^*\le\esssup X$, and in particular $x^*>-L$. Because $\Ip$ is convex, vanishes at $\Erw X$ and is nondecreasing, it is strictly increasing on $\{0<\Ip<\infty\}$; hence
\begin{equation}\label{eq:Ipstrict}
\Ip(y)<\ell\quad\text{for }y<x^*,\qquad \Ip(y)>\ell\quad\text{for }y>x^*.
\end{equation}

For $x<x^*$ and $t\in(0,1]$ let
\[
y_t(x):=\frac{x-(1-t)x^*}{t}=x^*-\frac{x^*-x}{t},
\]
and for $\rho\in(0,\infty]$ define (with $t\cdot\infty=\infty$ for $t>0$, so that $H(x;\infty)=\infty$)
\begin{equation}\label{eq:H}
H(x;\rho)=\inf_{t\in(0,1]}\Bigl\{t\rho+t\,\Imin\bigl(y_t(x)\bigr)\Bigr\}.
\end{equation}
The infimum in \eqref{eq:H} is unchanged if $\Imin$ is replaced by $I$: if $x<\Erw X$, then $y_t(x)<\Erw X$ for all $t$ and $I=\Imin$ there; if $x\ge\Erw X$, let $t_E$ be such that $y_{t_E}(x)=\Erw X$, and note that for $t>t_E$ one has $t\rho+tI(y_t)\ge t\rho\ge t_E\rho=t_E\rho+t_E\Imin(y_{t_E})$; with $\rho=-\log p_1$ it is therefore the rate function \eqref{eq:GH} of \cite{GantertHofelsauer2018}. We use $\Imin$ because it is monotone.

\begin{lemma}\label{lem:H}
Let $x\in(-L,x^*)$ and $\rho\in(0,\infty)$.
\begin{enumerate}
\item $0<H(x;\rho)\le\rho+\Imin(x)<\infty$, and the infimum in \eqref{eq:H} is attained at some $t_0\in(0,1]$.
\item $H(x;\cdot)$ is strictly increasing on $(0,\infty)$.
\item $H(\cdot;\rho)$ is convex and nonincreasing on $(-L,x^*)$; in particular it is continuous there.
\end{enumerate}
\end{lemma}

\begin{proof}
Let $f(t)=t\rho+t\Imin(y_t)$, $y_t=y_t(x)$. The map $t\mapsto y_t$ is continuous and increasing on $(0,1]$, with $y_1=x>-L$ and $y_t\to-\infty$ as $t\downarrow0$. If $L<\infty$, then $\Imin=\infty$ on $(-\infty,-L)$, so $f=\infty$ near $0$. If $L=\infty$, then for every $\lambda<0$, $t\Imin(y_t)\ge\lambda(tx^*-(x^*-x))-t\log\psi(\lambda)\to|\lambda|(x^*-x)$ as $t\downarrow0$, and $\lambda$ is arbitrary, so $f(t)\to\infty$. In both cases the infimum may be restricted to a compact interval $[t_1,1]\subset(0,1]$. On this interval $f$ is lower semicontinuous, so the infimum is attained at some $t_0$, and $f(t_0)\le f(1)=\rho+\Imin(x)<\infty$. Since $f(t_0)\ge t_0\rho>0$, (i) follows. For (ii), if $\rho<\rho'$ and $t_0'$ is a minimiser for $\rho'$, then $H(x;\rho)\le t_0'\rho+t_0'\Imin(y_{t_0'})<H(x;\rho')$. For (iii), $(t,z)\mapsto t\,\Imin(z/t)$ is the perspective of the convex function $\Imin$ and hence jointly convex on $(0,\infty)\times\R$; composing with the affine map $(t,x)\mapsto(t,x-(1-t)x^*)$ shows that $(t,x)\mapsto f$ is jointly convex, so $x\mapsto H(x;\rho)$ is convex. It is finite on $(-L,x^*)$ by (i), hence continuous there. Monotonicity follows from $y_t(x)$ being increasing in $x$ and $\Imin$ being nonincreasing.
\end{proof}

We say that the environment is of \emph{Schr\"oder type} if $p_1(\xi_1)>0$ a.s., equivalently $b(\xi_1)=1$ a.s. Otherwise $\PP_{\mathrm{env}}(b(\xi_1)\ge2)>0$, and we say that \emph{B\"ottcher-type generations occur}; generation $k$ is of Schr\"oder (B\"ottcher) type if $p_1(\xi_k)>0$ ($=0$). We call the environment \emph{mixed} if $0<\PP_{\mathrm{env}}(p_1(\xi_1)>0)<1$. Set
\[
\rhoq:=\E[-\log p_1(\xi_1)]\in(0,\infty],\qquad \rhoa:=-\log\E[p_1(\xi_1)]\in(0,\infty],\qquad \beta:=\E[\log b(\xi_1)]\in[0,\infty).
\]
Here $\rhoq>0$ and $\rhoa>0$ because $p_1(\xi_1)<1$ with positive probability, $\beta\le\ell<\infty$ because $b\le m$, and $\beta>0$ exactly when B\"ottcher-type generations occur. By Jensen's inequality $\rhoa\le\rhoq$, with strict inequality if $\rhoq<\infty$ and $p_1(\xi_1)$ is not a.s.\ constant.

\subsection{Main results}
The first-order behaviour of $M_n$ is classical in this setting (see \cite{MalleinMilos2019} and the references there). We record it because it follows from our estimates.

\begin{proposition}\label{prop:LLN}
Under Assumption~\ref{ass:standing}, for $\PP_{\mathrm{env}}$-a.e.\ $\xi$, $M_n/n\to x^*$ $P_\xi$-almost surely.
\end{proposition}

\begin{theorem}[Quenched upper bound and exponential regime]\label{thm:A}
Suppose Assumption~\ref{ass:standing} holds.
\begin{enumerate}
\item For $\PP_{\mathrm{env}}$-a.e.\ $\xi$ and every $x\in(-L,x^*)$,
\[
\limsup_{n\to\infty}\frac1n\log P_\xi(M_n\le xn)\le-H(x;\rhoq),
\]
where the right-hand side is $-\infty$ if $\rhoq=\infty$.
\item If $\rhoq<\infty$ (which implies that the environment is of Schr\"oder type), then for $\PP_{\mathrm{env}}$-a.e.\ $\xi$ and every $x\in(-L,x^*)$,
\[
\lim_{n\to\infty}\frac1n\log P_\xi(M_n\le xn)=-H(x;\rhoq).
\]
\end{enumerate}
\end{theorem}

\begin{theorem}[Quenched, B\"ottcher-type generations occur]\label{thm:B}
Suppose Assumption~\ref{ass:standing} holds, $\PP_{\mathrm{env}}(b(\xi_1)\ge2)>0$ and $L<\infty$. Then for $\PP_{\mathrm{env}}$-a.e.\ $\xi$ and every $x\in(-L,x^*)$,
\[
\lim_{n\to\infty}\frac1n\log\bigl[-\log P_\xi(M_n\le xn)\bigr]=\beta\,\tau^*(x),\qquad \tau^*(x):=\frac{x^*-x}{x^*+L}.
\]
\end{theorem}

For a constant environment $\nu=\delta_\zeta$ with $b(\zeta)=b\ge2$, Theorem~\ref{thm:B} recovers \eqref{eq:CH} under (A3). Together, Theorems~\ref{thm:A} and~\ref{thm:B} show that for a.e.\ $\xi$ the quenched probability decays at a finite exponential rate if and only if $\rhoq<\infty$, and doubly exponentially if B\"ottcher-type generations occur and $L<\infty$. For the annealed law we need, in part (ii), that the environment acts only through the shape of the offspring law and not through its mean.

\begin{theorem}[Annealed]\label{thm:annealed}
Suppose Assumption~\ref{ass:standing} holds and let $x\in(-L,x^*)$.
\begin{enumerate}
\item If $\PP_{\mathrm{env}}(p_1(\xi_1)>0)>0$, then
\[
\liminf_{n\to\infty}\frac1n\log\PP(M_n\le xn)\ge-H(x;\rhoa).
\]
\item Suppose (without the hypothesis of (i)) that there are constants $\bar m>1$ and $C_2<\infty$ such that, almost surely,
\begin{equation}\label{eq:constmean}
m(\xi_1)=\bar m,\qquad \sum_{k}k^2\,\xi_1(k)\le C_2 .
\end{equation}
Then $\limsup_{n}\frac1n\log\PP(M_n\le xn)\le-H(x;\rhoa)$ (which is $-\infty$ if $\rhoa=\infty$). In particular, if also $\PP_{\mathrm{env}}(p_1(\xi_1)>0)>0$, then
\[
\lim_{n\to\infty}\frac1n\log\PP(M_n\le xn)=-H(x;\rhoa).
\]
\end{enumerate}
\end{theorem}

In \eqref{eq:constmean} the non-random mean is essential for our proof of the upper bound, and we expect the annealed rate to be strictly smaller than $H(x;\rhoa)$ in general when the mean is random (see Section~\ref{sec:discussion}); the second-moment bound is used only to obtain $p_1\le\bar p<1$ and a uniform Paley--Zygmund estimate in Lemma~\ref{lem:uniform}. Condition \eqref{eq:constmean} implies $p_1(\xi_1)\le\bar p:=1-(\bar m-1)^2/C_2<1$: by the Cauchy--Schwarz inequality, $(\bar m-p_1)^2=(\sum_{k\ge2}k\,\xi_1(k))^2\le(1-p_1)C_2$ and $\bar m-p_1\ge\bar m-1$. The annealed probability never decays faster than exponentially in the regimes of Theorem~\ref{thm:annealed}(i), and it always decays at least exponentially:

\begin{proposition}\label{prop:annealedupper}
Under Assumption~\ref{ass:standing}, for every $x<x^*$ there is $c_x>0$ such that $\PP(M_n\le xn)\le e^{-c_xn}$ for all $n$ large.
\end{proposition}

\begin{corollary}[Quenched versus annealed]\label{cor:dichotomy}
Suppose Assumption~\ref{ass:standing} holds and let $x\in(-L,x^*)$.
\begin{enumerate}
\item If $\rhoq<\infty$ and $p_1(\xi_1)$ is not a.s.\ constant, then for a.e.\ $\xi$,
\[
\lim_{n\to\infty}\frac1n\log P_\xi(M_n\le xn)=-H(x;\rhoq)<-H(x;\rhoa)\le\liminf_{n\to\infty}\frac1n\log\PP(M_n\le xn).
\]
\item If $0<\PP_{\mathrm{env}}(p_1(\xi_1)>0)<1$ and $L<\infty$, then for a.e.\ $\xi$,
\[
P_\xi(M_n\le xn)=\exp\Bigl\{-e^{\,n\beta\tau^*(x)+o(n)}\Bigr\};
\]
moreover $e^{-nH(x;\rhoa)+o(n)}\le\PP(M_n\le xn)\le e^{-c_xn}$.
\end{enumerate}
\end{corollary}

\begin{proof}
(i) combines Theorem~\ref{thm:A}(ii), Theorem~\ref{thm:annealed}(i) and Lemma~\ref{lem:H}(ii) with $\rhoa<\rhoq$. (ii) combines Theorem~\ref{thm:B} (note $\PP_{\mathrm{env}}(b(\xi_1)\ge2)=\PP_{\mathrm{env}}(p_1(\xi_1)=0)>0$), Theorem~\ref{thm:annealed}(i) and Proposition~\ref{prop:annealedupper}.
\end{proof}

\begin{remark}\label{rem:interpretation}
In Corollary~\ref{cor:dichotomy}(ii) the annealed lower bound is obtained from an atypical environment: all of the first $t_0n$ generations ($t_0$ a minimiser in \eqref{eq:H} for $\rho=\rhoa$) are of Schr\"oder type, and, heuristically, the environment is tilted by $p_1$, which makes generations with large $p_1$ more likely. Under the quenched law the environment is typical; consequently a positive fraction $\PP_{\mathrm{env}}(b\ge2)$ of the first generations forces the population to grow at least like $e^{k\beta}$, and every one of these particles must then stay low. This is the mechanism behind the doubly exponential decay.
\end{remark}

\begin{remark}
Theorem~\ref{thm:A}(i) does not require the environment to be of Schr\"oder type. If $\rhoq=\infty$, which happens if B\"ottcher-type generations occur but also for some Schr\"oder-type environments with $\E[-\log p_1(\xi_1)]=\infty$, Theorem~\ref{thm:A}(i) shows that the decay is faster than exponential, also when the steps are unbounded from below. Theorem~\ref{thm:B} identifies the double-exponential order when B\"ottcher-type generations occur and $L<\infty$; the order of decay in the remaining cases is not determined in this paper.
\end{remark}

\section{Preliminary estimates}\label{sec:prelim}

Throughout, Assumption~\ref{ass:standing} is in force. We use $o(n)$ for quantities $a_n(\xi)$ with $a_n/n\to0$ for a.e.\ $\xi$, uniformly in the indices (shifts, lengths) indicated in each statement; constants denoted $c,C$ (with or without indices) do not depend on $n$ or $\xi$, but may depend on $\nu$, $\mu$ and the fixed parameters ($\eta,\eps,\delta,B,\dots$). In exponents, $e^{-\infty}=0$. Each lemma below that involves $\xi$ holds for a.e.\ $\xi$ for fixed parameters (Lemma~\ref{lem:harmonic} holds for every deterministic environment with $\omega_i(0)=0$, and Lemma~\ref{lem:uniform} for every $\omega\in\cS^{\N}$); in the proofs of the theorems the parameters range over fixed countable sequences, so the exceptional null sets can be combined into one.

\subsection{Ergodic sums over windows}

\begin{lemma}\label{lem:windows}
Let $(V_i)_{i\ge1}$ be i.i.d.\ real random variables with $\E|V_1|<\infty$, and let $C_0\ge1$. Then almost surely
\[
\lim_{n\to\infty}\frac1n\max_{0\le j<k\le C_0n}\Bigl|\sum_{i=j+1}^kV_i-(k-j)\E V_1\Bigr|=0,
\qquad
\lim_{n\to\infty}\frac1n\max_{1\le i\le C_0n}|V_i|=0 .
\]
\end{lemma}

\begin{proof}
Let $R_k=\sum_{i\le k}V_i-k\E V_1$. By the strong law, $R_k=o(k)$, hence $\max_{k\le C_0n}|R_k|=o(n)$, and the first claim follows from $|R_k-R_j|\le2\max_{l\le C_0n}|R_l|$. The second follows from $V_i=R_i-R_{i-1}+\E V_1$.
\end{proof}

We apply Lemma~\ref{lem:windows} to $\log m(\xi_i)$, $\log h_r(\xi_i)$ (see below), $\log p_1(\xi_i)$ when $\rhoq<\infty$, $\log b(\xi_i)$, and to block functionals of $\xi$.

\subsection{Random walk estimates}
For every $n\ge1$ and $a\in\R$, the Chernoff bounds give
\begin{equation}\label{eq:chernoff}
\Prw(S_n\le a)\le e^{-n\Imin(a/n)},\qquad \Prw(S_n\ge a)\le e^{-n\Ip(a/n)} .
\end{equation}
By Cram\'er's theorem \cite[Theorem~2.2.3]{DemboZeitouni2010}, for $y\in(-L,\infty)$,
\[
\liminf_{n\to\infty}\frac1n\log\Prw(S_n\le yn)\ge-\inf_{z<y}I(z)=-\Imin(y),
\]
where the last equality holds because $I=\Imin$ on $(-\infty,\Erw X]$, $I\ge0=\Imin$ on $[\Erw X,\infty)$, and $\Imin$ is continuous on the interior $(-L,\infty)$ of its domain; similarly $\liminf_n\frac1n\log\Prw(S_n\ge yn)\ge-\Ip(y)$ for $y<\esssup X$.

\subsection{Many-to-one and the first-born line}
Since the displacements are independent of the genealogy, for every measurable $g\ge0$,
\begin{equation}\label{eq:many}
E_\xi\Bigl[\sum_{u\in D_n}g(S_u)\Bigr]=\Bigl(\prod_{i=1}^nm(\xi_i)\Bigr)\Erw[g(S_n)] .
\end{equation}
By (A1) every particle has a child. Following the first-born child in every generation defines a particle $u^\circ_n\in D_n$ for each $n$; conditionally on the genealogical tree, $(S_{u^\circ_n})_n$ is a random walk with step law $\mu$. In particular, for any event $G$ measurable with respect to the genealogy, $P_\xi(G,\ S_{u^\circ_n}\le a)=P_\xi(G)\,\Prw(S_n\le a)$.

The upper tail of $M_n$ is controlled by \eqref{eq:many}, \eqref{eq:chernoff}, \eqref{eq:Ipstrict} and Lemma~\ref{lem:windows}:

\begin{lemma}\label{lem:uppertail}
For every $\eta>0$ and $\eps\in(0,1)$, for a.e.\ $\xi$,
\[
\max_{0\le j\le n,\ \eps n\le k\le n}P_{\theta^j\xi}\bigl(M_k>(x^*+\eta)k\bigr)\le e^{-\eps n(\Ip(x^*+\eta)-\ell)+o(n)},
\]
where $\Ip(x^*+\eta)-\ell\in(0,\infty]$; in particular the left-hand side tends to $0$.
\end{lemma}

\begin{proof}
By \eqref{eq:many} and \eqref{eq:chernoff}, $P_{\theta^j\xi}(M_k>(x^*+\eta)k)\le\exp\{\sum_{i=j+1}^{j+k}\log m(\xi_i)-k\Ip(x^*+\eta)\}$. By Lemma~\ref{lem:windows} the sum is $k\ell+o(n)$ uniformly, and $\Ip(x^*+\eta)>\ell$ by \eqref{eq:Ipstrict}. Since $k\ge\eps n$, the claim follows.
\end{proof}

\subsection{Harmonic moments}
For $r>1$ and $\zeta\in\cP$ set
\[
h_r(\zeta):=p_1(\zeta)+\bigl(1-p_1(\zeta)\bigr)2^{1-r}\in[2^{1-r},1].
\]

\begin{lemma}\label{lem:harmonic}
For every deterministic sequence $\omega=(\omega_1,\omega_2,\dots)\in\cP^{\N}$ with $\omega_i(0)=0$ for all $i$, every $r>1$ and $n\ge0$,
\[
E_\omega\bigl[Z_n^{-r}\bigr]\le\prod_{i=1}^nh_r(\omega_i),\qquad\text{hence}\qquad P_\omega(Z_n<K)\le K^r\prod_{i=1}^nh_r(\omega_i)\quad(K>0).
\]
\end{lemma}

\begin{proof}
Let $\phi_n(\omega)=E_\omega[Z_n^{-r}]$; $\phi_0=1$. Condition on the number $N_0$ of children of the root. If $N_0=1$, $Z_n$ has the law of $Z_{n-1}$ under $P_{\theta\omega}$. If $N_0\ge2$, then $Z_n\ge Z^{(1)}+Z^{(2)}\ge2\min(Z^{(1)},Z^{(2)})$, where $Z^{(1)},Z^{(2)}$ are the generation-$n$ sizes of the subtrees of the first two children, which are independent with the law of $Z_{n-1}$ under $P_{\theta\omega}$. Using $\min(a,b)^{-r}\le a^{-r}+b^{-r}$,
\[
\phi_n(\omega)\le p_1(\omega_1)\phi_{n-1}(\theta\omega)+\bigl(1-p_1(\omega_1)\bigr)2^{-r}\cdot2\,\phi_{n-1}(\theta\omega)=h_r(\omega_1)\phi_{n-1}(\theta\omega).
\]
Iterate, and use Markov's inequality $P_\omega(Z_n<K)\le K^rE_\omega[Z_n^{-r}]$.
\end{proof}

Let $\rho_r:=-\E[\log h_r(\xi_1)]$. Since $\log h_r\in[(1-r)\log2,0]$, Lemma~\ref{lem:windows} applies, and $\rho_r>0$ because $h_r(\xi_1)<1$ when $p_1(\xi_1)<1$. As $r\uparrow\infty$, $h_r\downarrow p_1$, so by monotone convergence
\begin{equation}\label{eq:rhor}
\rho_r\uparrow\rhoq,\qquad -\log\E[h_r(\xi_1)]\uparrow\rhoa .
\end{equation}

\subsection{Reaching the speed is not exponentially unlikely}
The next lemma is the main technical input. It holds uniformly over time-shifts of the environment.

\begin{lemma}\label{lem:chain}
For every $\eta\in(0,1]$ (the restriction is used only to simplify a constant in Step~3), $\delta>0$ and $\eps\in(0,1)$, for a.e.\ $\xi$ there is $n_0$ such that for all $n\ge n_0$,
\[
\min_{0\le j\le 2n,\ \eps n\le k\le n}P_{\theta^j\xi}\bigl(M_k\ge(x^*-\eta)k\bigr)\ge e^{-\delta n}.
\]
\end{lemma}

\begin{proof}
\emph{Step 1: an embedded branching process.} Let $y=x^*-\eta/4$. By \eqref{eq:Ipstrict}, $\Ip(y)<\ell$, and since $y<\esssup X$, Cram\'er's theorem yields $B'\in\N$ with
\begin{equation}\label{eq:Bprime}
\log\Prw(S_{B'}\ge yB')\ge-B'\bigl(\Ip(y)+\tfrac12(\ell-\Ip(y))\bigr).
\end{equation}
Given an environment, define $\mathcal Y_0=\{\varnothing\}$ (the root) and, recursively, $\mathcal Y_{g+1}$ as the set of $w\in D_{(g+1)B'}$ whose ancestor $v$ in generation $gB'$ belongs to $\mathcal Y_g$ and satisfies $S_w-S_v\ge yB'$. By the branching property, given $\xi$ the particles of $\mathcal Y_g$ produce their descendants in $\mathcal Y_{g+1}$ independently, so $Y_g=\#\mathcal Y_g$ is a branching process whose offspring law in generation $g$ depends only on the block $(\xi_{gB'+1},\dots,\xi_{(g+1)B'})$. Under $\PP_{\mathrm{env}}$ these blocks are i.i.d., so $(Y_g)$ is a branching process in i.i.d.\ random environment. By \eqref{eq:many}, its mean offspring number is $m^Y_g=\prod_{i=gB'+1}^{(g+1)B'}m(\xi_i)\,\Prw(S_{B'}\ge yB')$, and by \eqref{eq:Bprime},
\[
\E[\log m^Y_0]\ge B'\ell-B'\bigl(\Ip(y)+\tfrac12(\ell-\Ip(y))\bigr)=\tfrac12B'(\ell-\Ip(y))>0 ,
\]
while $\log m^Y_0\ge\log\Prw(S_{B'}\ge yB')>-\infty$ and $\E[\log m^Y_0]\le B'\ell<\infty$. Moreover, following the first-born line, $P_\xi(Y_1\ge1)\ge\Prw(S_{B'}\ge yB')>0$, so $\E[-\log P_\xi(Y_1\ge1)]<\infty$. We now apply \cite[Theorem~3]{AthreyaKarlin1971} (see also \cite{SmithWilkinson1969}): for a branching process in a stationary ergodic environment whose offspring generating functions $\varphi$ satisfy $\E[-\log(1-\varphi(0))]<\infty$ and $\E[(\log\varphi'(1))^-]<\E[(\log\varphi'(1))^+]\le\infty$, the quenched extinction probability is $<1$ almost surely. Here the environment variable of the generation-$g$ step of $Y$ is the block $(\xi_{gB'+1},\dots,\xi_{(g+1)B'})\in\cP^{B'}$, of which the offspring law of $Y$ is a measurable function; these blocks are i.i.d.\ under $\PP_{\mathrm{env}}$, the offspring law of $Y$ has finite mean since $Y_1\le Z_{B'}$, $\varphi'(1)=m^Y_0$ with $\log m^Y_0$ bounded below, $\E|\log m^Y_0|<\infty$ and $\E[\log m^Y_0]>0$ (so only the case $\E|\log\varphi'(1)|<\infty$ of the theorem is needed), and $1-\varphi(0)=P_\xi(Y_1\ge1)$, so both conditions hold. Hence the quenched extinction probability $q_Y(\xi)$ of $(Y_g)$ satisfies $q_Y(\xi)<1$ for a.e.\ $\xi$.

\emph{Step 2: block probabilities.} For $B\in\N$ and an environment $\omega$ let $a_B(\omega):=P_\omega(M_B\ge(x^*-\eta/2)B)$. Write $B=N'B'+R'$ with $0\le R'<B'$. On $\{Y_{N'}>0\}$ there is $w\in D_{N'B'}$ with $S_w\ge N'B'y$; following the first-born line of $w$ for $R'$ further generations,
\[
M_B\ge N'B'y+(S_{w'}-S_w),
\]
where $w'$ is the first-born descendant of $w$ in generation $B$ and, conditionally on $\cF_{N'B'}$, $S_{w'}-S_w$ has the law of $S_{R'}$. Let $c_1:=\min_{0\le R<B'}\Prw(S_R\ge R\,\Erw X-B')>0$. Since $P_\omega(Y_{N'}>0)\ge1-q_Y(\omega)$, and $w$ can be chosen $\cF_{N'B'}$-measurably on $\{Y_{N'}>0\}$, we get $M_B\ge yB-B'(|y|+|\Erw X|+1)\ge(x^*-\eta/2)B$ with probability at least $c_1(1-q_Y(\omega))$, as soon as $B\ge B_0:=4B'(|y|+|\Erw X|+1)/\eta$. Also, following a single line, $a_B(\omega)\ge\Prw(S_B\ge(x^*-\eta/2)B)\ge e^{-CB}$ for $B\ge B_0$ (after enlarging $B_0$), with $C:=\Ip(x^*-\eta/2)+1$. Hence, for $B\ge B_0$,
\[
0\le\frac{-\log a_B(\xi)}{B}\le\min\Bigl\{C,\ \frac{-\log c_1-\log(1-q_Y(\xi))}{B}\Bigr\}.
\]
Since $q_Y(\xi)<1$ a.s., the right-hand side tends to $0$ a.s.\ as $B\to\infty$ and is bounded by $C$; by dominated convergence, $\E[-\log a_B(\xi)]/B\to0$. Fix $B\ge B_0$ with $\E[-\log a_B(\xi)]\le\delta B/4$.

\emph{Step 3: chaining.} Fix $j\le2n$ and $k\in[\eps n,n]$, and write $k=NB+R$ with $0\le R<B$. Starting from the root of the subtree (at time $j$), choose recursively a particle $v_{i+1}$ at time $j+(i+1)B$ among the descendants of $v_i$ with $S_{v_{i+1}}-S_{v_i}\ge(x^*-\eta/2)B$, if such a particle exists. By the branching property, the probability that this succeeds for $i=0,\dots,N-1$ equals $\prod_{i=0}^{N-1}a_B(\theta^{j+iB}\xi)$. On success, following the first-born line of $v_N$ for $R$ generations gives, with conditional probability at least $c_2:=\min_{0\le R<B}\Prw(S_R\ge R\,\Erw X-B)>0$,
\[
M_k\ge N B(x^*-\eta/2)+R\,\Erw X-B\ge(x^*-\eta/2)k-B(|x^*|+|\Erw X|+2)\ge(x^*-\eta)k
\]
for $n$ large (as $k\ge\eps n$). Therefore
\[
P_{\theta^j\xi}\bigl(M_k\ge(x^*-\eta)k\bigr)\ge c_2\exp\Bigl\{-\sum_{i=0}^{N-1}V(\theta^{j+iB}\xi)\Bigr\},\qquad V(\omega):=-\log a_B(\omega).
\]

\emph{Step 4: uniformity.} Write $j=j_0+i_0B$ with $0\le j_0<B$. For fixed $j_0$, the variables $V^{(j_0)}_i:=V(\theta^{j_0+iB}\xi)$, $i\ge0$, depend on the disjoint blocks $(\xi_{j_0+iB+1},\dots,\xi_{j_0+(i+1)B})$ and are therefore i.i.d., nonnegative, bounded by $CB$, with mean at most $\delta B/4$. The indices involved satisfy $i\le(j+k)/B\le3n/B$, so Lemma~\ref{lem:windows} (with $C_0=3$), applied to each of the $B$ sequences $(V^{(j_0)}_i)_i$, gives, uniformly in $j\le2n$ and $k\le n$,
\[
\sum_{i=0}^{N-1}V(\theta^{j+iB}\xi)=\sum_{i=i_0}^{i_0+N-1}V^{(j_0)}_i\le N\frac{\delta B}{4}+o(n)\le\frac{\delta n}{4}+o(n).
\]
Hence the probability is at least $c_2e^{-\delta n/4-o(n)}\ge e^{-\delta n}$ for $n$ large.
\end{proof}

\subsection{An a priori exponential bound}

\begin{lemma}\label{lem:apriori}
For every $\eta>0$ and $\eps\in(0,1)$ there is $c=c(\eta,\eps)>0$ such that for a.e.\ $\xi$ and all large $n$,
\[
\max_{0\le j\le n,\ \eps n\le k\le n}P_{\theta^j\xi}\bigl(M_k\le(x^*-\eta)k\bigr)\le e^{-cn}.
\]
\end{lemma}

\begin{proof}
We may assume $\eta\le1$. Choose $A>0$ as follows: if $L<\infty$, $A:=\max(L,0)+1$; otherwise $A$ with $\Imin(-A)\ge\ell+1$ (possible since $\Imin(y)\to\infty$ as $y\to-\infty$). Choose $\sigma\in(0,\frac12)$ with $\sigma(A+|x^*|+1)\le\eta/4$, put $\kappa:=\sigma\eps\rho_2/4$, and fix $j\le n$, $k\in[\eps n,n]$ and $s:=\lfloor\sigma k\rfloor$. Writing $P_j=P_{\theta^j\xi}$,
\[
\{M_k\le(x^*-\eta)k\}\subset E_1\cup E_2\cup E_3,
\]
where $E_1=\{Z_s<e^{\kappa n}\}$, $E_2=\{\min_{u\in D_s}S_u<-As\}$ and
\[
E_3=\bigl\{Z_s\ge e^{\kappa n},\ \min_{u\in D_s}S_u\ge-As,\ M^u_{k-s}\le(x^*-\eta)k+As\ \ \forall u\in D_s\bigr\}.
\]
By Lemma~\ref{lem:harmonic} with $r=2$ and Lemma~\ref{lem:windows}, $P_j(E_1)\le e^{2\kappa n}\exp\{-s\rho_2+o(n)\}\le\exp\{2\kappa n-4\kappa n+o(n)\}$, because $s\ge\sigma\eps n-1$. If $L<\infty$ then $P_j(E_2)=0$; otherwise, by \eqref{eq:many}, \eqref{eq:chernoff} and Lemma~\ref{lem:windows}, $P_j(E_2)\le\exp\{s\ell+o(n)-s\Imin(-A)\}\le\exp\{-s+o(n)\}$. Finally,
\[
(x^*-\eta/2)(k-s)-\bigl[(x^*-\eta)k+As\bigr]=\tfrac{\eta}{2}k-s(x^*-\tfrac{\eta}{2}+A)\ge\tfrac{\eta}{2}k-\sigma k(|x^*|+1+A)\ge\tfrac{\eta}{4}k>0,
\]
so on $E_3$ each of the at least $e^{\kappa n}$ subtrees rooted in $D_s$ fails to reach $(x^*-\eta/2)(k-s)$ within $k-s$ generations. Conditionally on $\cF_s$ these events are independent, and by Lemma~\ref{lem:chain} (applied with $\eta/2$, $\delta=\kappa/2$, $\eps/2$ and the shift $j+s\le2n$; note $k-s\ge k/2\ge\eps n/2$), each has probability at most $1-e^{-\kappa n/2}$. Hence $P_j(E_3)\le(1-e^{-\kappa n/2})^{e^{\kappa n}}\le\exp\{-e^{\kappa n/2}\}$. All bounds are uniform in $j,k$, and the claim follows with $c=\min\{\kappa,\sigma\eps/2\}$.
\end{proof}

\begin{proof}[Proof of Proposition~\ref{prop:LLN}]
By Lemma~\ref{lem:apriori} and Lemma~\ref{lem:uppertail} (with $j=0$, $k=n$), both $P_\xi(M_n\le(x^*-\eta)n)$ and $P_\xi(M_n>(x^*+\eta)n)$ are summable in $n$. Borel--Cantelli and $\eta\downarrow0$ along a sequence give the claim.
\end{proof}

\section{Proof of Theorem~\ref{thm:A}}\label{sec:proofA}

\subsection{Upper bound}
Fix $x\in(-L,x^*)$, $r>1$, $\eps\in(0,1)$, $\eta\in(0,(x^*-x)/2)$ and an integer $K\ge2$, and let $c=c(\eta,\eps)$ be given by Lemma~\ref{lem:apriori}. Put $n_\eps=\lfloor(1-\eps)n\rfloor$ and $T_K:=\inf\{k\ge0:Z_k\ge K\}$; since $Z_0=1<K$, $T_K\ge1$. Then
\[
P_\xi(M_n\le xn)\le P_\xi(T_K>n_\eps,\ M_n\le xn)+\sum_{k=1}^{n_\eps}P_\xi(T_K=k,\ M_n\le xn).
\]

\emph{Late growth.} Since $Z$ is nondecreasing by (A1), $\{T_K>n_\eps\}=\{Z_{n_\eps}<K\}$ is measurable with respect to the genealogy and $M_n\ge S_{u^\circ_n}$, Lemma~\ref{lem:harmonic} and \eqref{eq:chernoff} give
\[
P_\xi(T_K>n_\eps,\ M_n\le xn)\le P_\xi(Z_{n_\eps}<K)\,\Prw(S_n\le xn)\le K^r\prod_{i=1}^{n_\eps}h_r(\xi_i)\,e^{-n\Imin(x)} .
\]

\emph{Growth at time $k\le n_\eps$.} Let $a_k:=xn-(x^*-\eta)(n-k)$ and $g_k(z):=P_{\theta^k\xi}(M_{n-k}\le z)$. On $\{T_K=k\}\in\cF_k$, the branching property gives $P_\xi(M_n\le xn\mid\cF_k)=\prod_{u\in D_k}g_k(xn-S_u)$. If $S_u\ge a_k$, then $xn-S_u\le(x^*-\eta)(n-k)$ and, as $n-k\ge\eps n$, Lemma~\ref{lem:apriori} gives $g_k(xn-S_u)\le e^{-cn}$ for $n$ large, uniformly in $k$. Since $Z_k\ge K$ on $\{T_K=k\}$,
\[
\mathbf 1_{\{T_K=k\}}P_\xi(M_n\le xn\mid\cF_k)\le\mathbf 1_{\{Z_{k-1}<K,\ \exists u\in D_k:\,S_u<a_k\}}+e^{-cKn}.
\]
Taking expectations and using that the displacements are independent of the genealogy,
\begin{align*}
P_\xi(Z_{k-1}<K,\ \exists u\in D_k:S_u<a_k)&\le E_\xi\bigl[\mathbf 1_{\{Z_{k-1}<K\}}Z_k\bigr]\,\Prw(S_k<a_k)\\
&\le K\,m(\xi_k)\,P_\xi(Z_{k-1}<K)\,e^{-k\Imin(a_k/k)} .
\end{align*}
With Lemma~\ref{lem:harmonic}, we obtain for $1\le k\le n_\eps$
\[
P_\xi(T_K=k,\ M_n\le xn)\le K^{r+1}m(\xi_k)\prod_{i=1}^{k-1}h_r(\xi_i)\,e^{-k\Imin(a_k/k)}+e^{-cKn}.
\]

\emph{Collecting terms.} By Lemma~\ref{lem:windows}, $\sum_{i<k}\log h_r(\xi_i)=-k\rho_r+o(n)$ and $\log m(\xi_k)=o(n)$ uniformly in $k\le n$. Writing $t=k/n$ and noting $a_k/k=y^\eta_t:=\bigl(x-(1-t)(x^*-\eta)\bigr)/t$, we arrive at
\[
\limsup_{n\to\infty}\frac1n\log P_\xi(M_n\le xn)\le-\min\Bigl\{(1-\eps)\rho_r+\Imin(x),\ \inf_{t\in(0,1]}f_{r,\eta}(t),\ cK\Bigr\},
\]
where $f_{r,\eta}(t):=t\rho_r+t\Imin(y^\eta_t)$.
Since $c$ does not depend on $K$, we may let $K\to\infty$, and then $\eps\downarrow0$. As $f_{r,\eta}(1)=\rho_r+\Imin(x)$, this yields $\limsup_n\frac1n\log P_\xi(M_n\le xn)\le-\inf_{t\in(0,1]}f_{r,\eta}(t)$ for every $r>1$ and small $\eta>0$.

\emph{Removing $r$ and $\eta$.} Let $f(t)=t\rhoq+t\Imin(y_t(x))$. As $r\uparrow\infty$ and $\eta\downarrow0$, $f_{r,\eta}(t)\uparrow f(t)$ for every $t\in(0,1]$: indeed $\rho_r\uparrow\rhoq$ by \eqref{eq:rhor}, and $y^\eta_t\downarrow y_t(x)$ while $\Imin$ is nonincreasing and lower semicontinuous, so $\Imin(y^\eta_t)\uparrow\Imin(y_t(x))$. Fix $\Gamma<\infty$. Since $y^\eta_t\le x^*-(x^*-x)/(2t)$ for $\eta<(x^*-x)/2$ and $f_{r,\eta}(t)\ge t\Imin(y^\eta_t)$, the argument in the proof of Lemma~\ref{lem:H}(i) gives $t_1=t_1(\Gamma)>0$ such that $f_{r,\eta}(t)\ge\Gamma$ for all $t\in(0,t_1]$, all $r>1$ and all $\eta<(x^*-x)/2$. On the compact interval $[t_1,1]$, the functions $f_{r,\eta}$ are lower semicontinuous and increase to $f$ along any sequence $r_i\uparrow\infty$, $\eta_i\downarrow0$. Hence $\inf_{[t_1,1]}f_{r_i,\eta_i}\to\inf_{[t_1,1]}f$: otherwise there are $c'<\inf_{[t_1,1]}f$ and $t_i\in[t_1,1]$ with $f_{r_i,\eta_i}(t_i)\le c'$; along a subsequence $t_i\to\bar t$, and for fixed $i_0$ and $i\ge i_0$, $f_{r_{i_0},\eta_{i_0}}(t_i)\le f_{r_i,\eta_i}(t_i)\le c'$, so lower semicontinuity gives $f_{r_{i_0},\eta_{i_0}}(\bar t)\le c'$ and, letting $i_0\to\infty$, $f(\bar t)\le c'$, a contradiction. Consequently $\lim_i\inf_{(0,1]}f_{r_i,\eta_i}\ge\min\{\Gamma,\inf_{[t_1,1]}f\}\ge\min\{\Gamma,H(x;\rhoq)\}$, and since $\Gamma$ is arbitrary the limit is at least $H(x;\rhoq)$ (and is infinite if $\rhoq=\infty$). This proves the upper bound in (i) for fixed $x$; since $P_\xi(M_n\le xn)$ is nondecreasing in $x$ and $H(\cdot;\rhoq)$ is continuous (Lemma~\ref{lem:H}(iii)), a countable dense set of $x$ suffices to obtain the statement for all $x$ simultaneously. If $\rhoq=\infty$, the bound $-\infty$ holds for all rational $x<x^*$ simultaneously, and for general $x$ it follows from $P_\xi(M_n\le xn)\le P_\xi(M_n\le x'n)$ with rational $x'\in(x,x^*)$.

\subsection{Lower bound}
Assume now $\rhoq<\infty$, fix $x\in(-L,x^*)$, $t\in(0,1)$ with $y_t(x)>-L$, and $\eta>0$ small enough that $y^{-\eta}_t:=\bigl(x-(1-t)(x^*+\eta)\bigr)/t>-L$. Let $k=\lfloor tn\rfloor$ and $z_n:=xn-(x^*+\eta)(n-k)$. Consider the event that $Z_k=1$, that the unique particle $v$ of generation $k$ satisfies $S_v\le z_n$, and that $M^v_{n-k}\le(x^*+\eta)(n-k)$. On this event $M_n\le xn$. Since $\{Z_k=1\}$ depends only on the genealogy, $S_v=S_{u^\circ_k}$ on it, and the subtree of $v$ is independent of $\cF_k$,
\[
P_\xi(M_n\le xn)\ge\prod_{i=1}^kp_1(\xi_i)\cdot\Prw(S_k\le z_n)\cdot\bigl(1-P_{\theta^k\xi}(M_{n-k}>(x^*+\eta)(n-k))\bigr).
\]
The last factor tends to $1$ by Lemma~\ref{lem:uppertail}. By the strong law, $\prod_{i\le k}p_1(\xi_i)=e^{-tn\rhoq+o(n)}$. Since $z_n/k\to y^{-\eta}_t>-L$, Cram\'er's theorem gives $\liminf_n\frac1n\log\Prw(S_k\le z_n)\ge-t\Imin(y^{-\eta}_{t}-\gamma)$ for every $\gamma\in(0,y^{-\eta}_t+L)$, hence $\ge-t\Imin(y^{-\eta}_t)$ by continuity of $\Imin$ on $(-L,\infty)$. Letting $\eta\downarrow0$ and using again the continuity of $\Imin$,
\[
\liminf_{n\to\infty}\frac1n\log P_\xi(M_n\le xn)\ge-\bigl(t\rhoq+t\Imin(y_t(x))\bigr).
\]
The case $t=1$ follows from the same argument with $k=n$ and no last factor. Since $t\mapsto t\rhoq+t\Imin(y_t(x))$ is continuous on $\{t:y_t(x)>-L\}$, it suffices to take rational $t$ and $\eta$, so that only countably many null sets are excluded. It remains to see that the supremum over $t$ with $y_t(x)>-L$ equals $-H(x;\rhoq)$. Let $t_0$ be a minimiser in \eqref{eq:H} (Lemma~\ref{lem:H}(i)). If $y_{t_0}(x)>-L$ there is nothing to show. Otherwise $y_{t_0}(x)=-L$ and $\Imin(-L)<\infty$; for $t\downarrow t_0$ we have $y_t(x)\downarrow-L$ with $y_t(x)>-L$, and writing $y_t(x)=(1-s)(-L)+s\,x$ with $s\downarrow0$, convexity gives $\Imin(y_t(x))\le(1-s)\Imin(-L)+s\Imin(x)$, so $\limsup_{t\downarrow t_0}\bigl(t\rhoq+t\Imin(y_t(x))\bigr)\le t_0\rhoq+t_0\Imin(-L)=H(x;\rhoq)$. The extension to all $x$ simultaneously is as before.

\section{Proof of Theorem~\ref{thm:B}}\label{sec:proofB}
Write $b_i=b(\xi_i)$ and $\Pi_k=\prod_{i=1}^kb_i$. Since every particle of generation $i-1$ has at least $b_i$ children, $Z_k\ge\Pi_k$ surely, and since $X\ge-L$ a.s., $S_u\ge-Lk$ for every $u\in D_k$. By the strong law, $\log\Pi_k=k\beta+o(k)$.

\subsection{Upper bound on the probability}
Fix $x\in(-L,x^*)$, $\eta\in(0,x^*-x)$, $\delta>0$, let $\tau_\eta:=\frac{x^*-\eta-x}{x^*-\eta+L}\in(0,1)$ and $k:=\lfloor\tau_\eta n\rfloor$. For every $u\in D_k$,
\[
xn-S_u\le xn+Lk\le(x^*-\eta)(n-k),
\]
where the second inequality is equivalent to $k\le\tau_\eta n$. Hence, by the branching property and Lemma~\ref{lem:chain} (applied with $\eta/2$, where we may assume $\eta\le2$, and $\eps=1-\tau_\eta$), writing $g_k(z)=P_{\theta^k\xi}(M_{n-k}\le z)$,
\[
P_\xi(M_n\le xn)=E_\xi\Bigl[\prod_{u\in D_k}g_k(xn-S_u)\Bigr]\le\Bigl(1-P_{\theta^k\xi}\bigl(M_{n-k}\ge(x^*-\tfrac\eta2)(n-k)\bigr)\Bigr)^{\Pi_k}\le\bigl(1-e^{-\delta n}\bigr)^{\Pi_k}
\]
for $n$ large. Thus $-\log P_\xi(M_n\le xn)\ge\Pi_ke^{-\delta n}$ and
\[
\liminf_{n\to\infty}\frac1n\log\bigl[-\log P_\xi(M_n\le xn)\bigr]\ge\tau_\eta\beta-\delta .
\]
Let $\delta\downarrow0$ and $\eta\downarrow0$.

\subsection{Lower bound on the probability}
For an integer $J\ge2$ and an offspring law $\zeta\in\cP$ with $\zeta(0)=0$ and $m(\zeta)<\infty$, let $d_J(\zeta):=\min\{d\ge1:\zeta([1,d])\ge1/J\}$. By Markov's inequality, $\zeta([1,d])\ge1-m(\zeta)/(d+1)\ge\frac12$ for $d=\lceil2m(\zeta)\rceil$, so $b(\zeta)\le d_J(\zeta)\le\lceil2m(\zeta)\rceil$; in particular $\log d_J(\xi_1)$ is integrable by (A2). Moreover $d_J(\zeta)$ is nonincreasing in $J$ and equals $b(\zeta)$ for $J$ large, so by dominated convergence
\begin{equation}\label{eq:betaJ}
\beta_J:=\E[\log d_J(\xi_1)]\downarrow\beta\qquad(J\to\infty).
\end{equation}
Fix $J\ge2$, $x\in(-L,x^*)$, $\eta>0$ and $\eps'\in(0,x+L)$, let $\varpi:=\Prw(X\le-L+\eps')>0$, $\Pi^J_k:=\prod_{i=1}^kd_J(\xi_i)$, and
\[
\tau':=\frac{x^*+\eta-x}{x^*+\eta+L-\eps'}\in(0,1),\qquad k:=\lceil\tau'n\rceil .
\]
Let $G$ be the event that for each $i\le k$, every particle of generation $i-1$ has at most $d_J(\xi_i)$ children and each of these children is displaced by at most $-L+\eps'$. On $G$, $Z_i\le \Pi^J_i$ for $i\le k$ and $S_u\le-(L-\eps')k$ for all $u\in D_k$. A particle of generation $i-1$ satisfies the requirement with conditional probability at least $\xi_i([1,d_J(\xi_i)])\,\varpi^{d_J(\xi_i)}\ge J^{-1}\varpi^{d_J(\xi_i)}$, independently of the others, and on $G$ there are at most $\Pi^J_{i-1}$ of them. Hence
\[
-\log P_\xi(G)\le\sum_{i=1}^k\Bigl(\Pi^J_{i-1}\log J+\Pi^J_i\log\tfrac1\varpi\Bigr)\le C_J\,k\,\Pi^J_k,\qquad C_J:=\log J+\log\tfrac1\varpi.
\]
On $G$, for every $u\in D_k$, $xn-S_u\ge xn+(L-\eps')k\ge(x^*+\eta)(n-k)$ by the choice of $\tau'$. Hence, by the branching property and Lemma~\ref{lem:uppertail},
\[
P_\xi(M_n\le xn)\ge P_\xi(G)\Bigl(1-P_{\theta^k\xi}\bigl(M_{n-k}>(x^*+\eta)(n-k)\bigr)\Bigr)^{\Pi^J_k}\ge P_\xi(G)\,2^{-\Pi^J_k}
\]
for $n$ large. Therefore $-\log P_\xi(M_n\le xn)\le\Pi^J_k\bigl(C_Jk+\log2\bigr)$, and since $\log\Pi^J_k=k\beta_J+o(n)$ by the strong law,
\[
\limsup_{n\to\infty}\frac1n\log\bigl[-\log P_\xi(M_n\le xn)\bigr]\le\tau'\beta_J .
\]
Letting $J\to\infty$ (using \eqref{eq:betaJ}) and then $\eta,\eps'\downarrow0$ gives the upper bound $\beta(x^*-x)/(x^*+L)$, since $\tau'\to(x^*-x)/(x^*+L)$. As in Section~\ref{sec:proofA}, monotonicity in $x$ and continuity of the limit yield the statement for all $x$ simultaneously.

\section{Proof of Theorem~\ref{thm:annealed}}\label{sec:proofannealed}

\subsection{Part (i)}
Let $t$, $\eta$, $k$ and $z_n$ be as in the lower bound of Section~\ref{sec:proofA}. The quenched bound obtained there holds for every environment:
\[
P_\xi(M_n\le xn)\ge\prod_{i=1}^kp_1(\xi_i)\cdot\Prw(S_k\le z_n)\cdot\bigl(1-P_{\theta^k\xi}(M_{n-k}>(x^*+\eta)(n-k))\bigr).
\]
The first factor depends only on $(\xi_1,\dots,\xi_k)$ and the last only on $(\xi_{k+1},\dots,\xi_n)$, so they are independent under $\PP_{\mathrm{env}}$. Hence
\[
\PP(M_n\le xn)\ge\bigl(\E[p_1(\xi_1)]\bigr)^k\,\Prw(S_k\le z_n)\,\E\bigl[1-P_{\theta^k\xi}(M_{n-k}>(x^*+\eta)(n-k))\bigr].
\]
The last expectation tends to $1$ by Lemma~\ref{lem:uppertail} and bounded convergence, and $(\E p_1(\xi_1))^k=e^{-tn\rhoa+O(1)}$. The rest of the argument is identical to Section~\ref{sec:proofA}.

\subsection{Part (ii)}
Let $\cS$ be the set of $\zeta\in\cP$ with $\zeta(0)=0$, $m(\zeta)=\bar m$ and $\sum_kk^2\zeta(k)\le C_2$ (hence $p_1(\zeta)\le\bar p<1$); by \eqref{eq:constmean}, $\xi_i\in\cS$ for all $i$ almost surely. Here $\ell=\log\bar m$. We first show that Lemmas~\ref{lem:chain} and \ref{lem:apriori} hold uniformly over \emph{all} deterministic environments in $\cS^{\N}$.

\begin{lemma}\label{lem:uniform}
Under \eqref{eq:constmean}:
\begin{enumerate}
\item for every $\eta>0$ there are $c_\eta>0$ and $k_0$ such that $P_\omega(M_k\ge(x^*-\eta)k)\ge c_\eta$ for all $\omega\in\cS^{\N}$ and $k\ge k_0$;
\item for every $\eta>0$ there are $c>0$ and $k_1$ such that $P_\omega(M_k\le(x^*-\eta)k)\le e^{-ck}$ for all $\omega\in\cS^{\N}$ and $k\ge k_1$.
\end{enumerate}
\end{lemma}

\begin{proof}
(i) Let $y=x^*-\eta/4$ and choose $B'$ such that $\mu_Y:=\bar m^{B'}\Prw(S_{B'}\ge yB')\ge2$, which is possible by \eqref{eq:Ipstrict} and Cram\'er's theorem. Define $(Y_g)$ as in Step~1 of the proof of Lemma~\ref{lem:chain}. Under $P_\omega$ it is a Galton--Watson process in the varying environment determined by $\omega$, with mean $\mu_Y$ in every generation and offspring variances $s_g^2\le E_{\theta^{gB'}\omega}[Z_{B'}^2]\le C_3$, where $C_3$ depends only on $B'$, $\bar m$ and $C_2$ (iterate $E_\omega[Z_i^2]\le\bar m^2E_\omega[Z_{i-1}^2]+C_2E_\omega[Z_{i-1}]$). With $W_g=Y_g/\mu_Y^g$ one has $E_\omega[W_g^2]=E_\omega[W_{g-1}^2]+s_{g-1}^2\mu_Y^{-(g+1)}$, hence $E_\omega[W_g^2]\le1+C_3$ and, by the Paley--Zygmund inequality, $P_\omega(Y_g>0)\ge(E_\omega W_g)^2/E_\omega[W_g^2]\ge(1+C_3)^{-1}$ for all $g$. The remainder of $k$ modulo $B'$ is handled by the first-born line as in Step~2 of Lemma~\ref{lem:chain}.

(ii) Repeat the proof of Lemma~\ref{lem:apriori} with $k$ in place of $n$ and no shift (so that $k$ plays the role of both $n$ and the window length), with $A$ and $\sigma$ chosen as there, but with $\kappa:=\sigma\log(1/\bar h)/4>0$, where $\bar h:=\bar p+(1-\bar p)/2<1$. By Lemma~\ref{lem:harmonic}, $P_\omega(Z_s<e^{\kappa k})\le e^{2\kappa k}\bar h^{\,s}\le\bar h^{-1}e^{-2\kappa k}$ since $s\ge\sigma k-1$; by \eqref{eq:many}, $P_\omega(E_2)\le\bar m^s\Prw(S_s<-As)\le e^{-s}$; and by (i), $P_\omega(E_3)\le(1-c_{\eta/2})^{e^{\kappa k}}$; the margin $\eta k/4$ in the definition of $E_3$ is as before. All bounds are uniform in $\omega\in\cS^{\N}$, and the claim holds with $c:=\min\{\kappa,\sigma/2\}$.
\end{proof}

Fix an environment $\omega\in\cS^{\N}$ and repeat the upper bound of Section~\ref{sec:proofA}, using Lemma~\ref{lem:uniform}(ii) in place of Lemma~\ref{lem:apriori}. Since Lemma~\ref{lem:uniform}(ii) is applied to windows of length $n-k\ge\eps n$, the term $e^{-cKn}$ of Section~\ref{sec:proofA} becomes $e^{-c\eps Kn}$. This yields, for $n\ge n_1:=\lceil k_1/\eps\rceil$ (independent of $\omega$),
\[
P_\omega(M_n\le xn)\le K^r\prod_{i=1}^{n_\eps}h_r(\omega_i)\,e^{-n\Imin(x)}+\sum_{k=1}^{n_\eps}\Bigl(K^{r+1}\bar m\prod_{i=1}^{k-1}h_r(\omega_i)\,e^{-k\Imin(a_k/k)}+e^{-c\eps Kn}\Bigr).
\]
Integrating over $\omega=\xi$ and using independence, $\E\prod_{i<k}h_r(\xi_i)=(\E h_r(\xi_1))^{k-1}=e^{-(k-1)\rho_{\mathrm a,r}}$ with $\rho_{\mathrm a,r}:=-\log\E[h_r(\xi_1)]$. The optimisation of Section~\ref{sec:proofA} applies verbatim with $\rho_r$ replaced by $\rho_{\mathrm a,r}$, and $\rho_{\mathrm a,r}\uparrow\rhoa$ by \eqref{eq:rhor}. This gives $\limsup_n\frac1n\log\PP(M_n\le xn)\le-H(x;\rhoa)$, and together with part (i) the limit.

\subsection{Proof of Proposition~\ref{prop:annealedupper}}
By Proposition~\ref{prop:LLN}, $M_B/B\to x^*$ $\PP$-almost surely as $B\to\infty$. Since $M_B\ge S_{u^\circ_B}$ and, by Jensen's inequality, $e^{-2S_B/B}\le B^{-1}\sum_{i\le B}e^{-2X_i}$,
\[
\sup_B\E_{\PP}\bigl[e^{-2M_B/B}\bigr]\le\Erw\bigl[e^{-2X}\bigr]<\infty
\]
by (A3). Hence $(e^{-M_B/B})_B$ is uniformly integrable and $\E_\PP[e^{x-M_B/B}]\to e^{x-x^*}<1$. Fix $B$ with $\varrho_B:=\E_\PP[e^{x-M_B/B}]<1$.

Starting from the root, choose recursively $v_{i+1}$ as a rightmost descendant at time $(i+1)B$ of $v_i$, and put $\Delta_i:=S_{v_{i+1}}-S_{v_i}$. Given $\xi$ and $\cF_{iB}$, $\Delta_i$ has the law of $M_B$ under $P_{\theta^{iB}\xi}$, which depends only on $(\xi_{iB+1},\dots,\xi_{(i+1)B})$; since these blocks are i.i.d., $\E_\PP[\prod_{i<N}e^{-\Delta_i/B}]=(\E_\PP e^{-M_B/B})^N$. For $n=NB+R$ with $0\le R<B$, following the first-born line of $v_N$ for $R$ generations gives $M_n\ge\sum_{i<N}\Delta_i+S'_R$, where $S'_R$ has the law of $S_R$ and is independent of $(\Delta_i)_{i<N}$. By the exponential Markov inequality,
\[
\PP(M_n\le xn)\le e^{xn/B}\,\E_\PP\Bigl[\prod_{i<N}e^{-\Delta_i/B}\Bigr]\Erw[e^{-S_R/B}]\le \varrho_B^{\,N}\,e^{|x|}\max_{R<B}\Erw[e^{-S_R/B}],
\]
which decays exponentially in $n$.

\section{Numerical illustration}\label{sec:numerics}

\subsection{Computational methods}
The scripts, the implementation of the recursion below and the numerical data are provided in Online Resource~1.

\emph{Algorithm.} For integer-valued steps, the law of $M_n$ under $P_\omega$ can be computed without simulation. Let $u_k(y)=P_{\theta^k\omega}(M_{n-k}\le y)$ and let $f_k$ be the generating function of $\omega_k$. The branching property gives the backward recursion
\[
u_n(y)=\mathbf 1_{\{y\ge0\}},\qquad u_k(y)=f_{k+1}\Bigl(\sum_{z}\mu(\{z\})\,u_{k+1}(y-z)\Bigr),\qquad P_\omega(M_n\le y)=u_0(y).
\]
Since $p_0=0$ and the steps take values in $\{-1,0,1\}$ in our examples, $u_k(y)=1$ for $y\ge n-k$ and $u_k(y)=0$ for $y<-(n-k)$, so the recursion on a finite grid involves no truncation. We carry both $-\log u_k$ and $\log(1-u_k)$ in double precision, so that doubly exponentially small probabilities and tail probabilities far below $10^{-308}$ remain representable. Carrying $\log(1-u_k)$ separately is necessary: rounding errors in $1-u_k$ are amplified by the mean offspring number in each generation.

\emph{Validation.} We validated the implementation against a $40$-digit reference computation that uses exact rational inputs and carries $1-u$ separately (rebuilding $u$ from $1-u$ where $1-u<\frac12$, for the reason just given). For the first $200$ generations of the seed-$11$ environments of examples (A) and (B) below, and for the constant environment $\zeta(1)=\zeta(3)=\frac12$ with $n=1600$, the errors of $-\log u_0$ and $\log(1-u_0)$ on the whole grid are at most $3.7\cdot10^{-14}\max(1,|\text{value}|)$; the smallest tail probability in these comparisons is about $10^{-282}$. As a sanity check, a Monte Carlo simulation with $2\cdot10^4$ samples of $M_6$ in the seed-$11$ environments of (A), (B) and (C) agrees with the computed distribution function at all $12$ grid points with $0<P<1$ in each of the three environments (largest $|z|$-score $1.53$).

\emph{Constant-environment benchmark.} In the constant environment $\zeta(1)=\zeta(3)=\frac12$ with $y=\lfloor0.3n\rfloor$, the computed values $w_n=-\log P(M_n\le y)$ for $n=800,1600,3200,6400$ fit $w_n=\hat Hn+\hat a\log n+\hat c$ with $\hat H=0.366063$ and $\hat a\approx-0.42$, against the Gantert--H\"ofelsauer rate $H(0.3)=0.366066$. This fit has one residual degree of freedom, and all four three-point fits give $\hat H\in[0.366061,0.366063]$, slightly below the theoretical value: $\hat H-H(0.3)\approx-3\cdot10^{-6}$, presumably from corrections of lower order than $\log n$ that the three-parameter model does not represent and that are partly absorbed into $\hat H$ and $\hat a$; the spread of the three-point fits is therefore not an error bar, and $\hat a$ is an empirical value that depends on which lower-order terms are included in the fit.

\emph{Annealed law.} The annealed probability is computed for $n\le16$ by summing over all $2^n$ environment sequences. For larger $n$ we use importance sampling over environments: for $x=0$ and $k$ the integer nearest to $t_0n$, where $t_0$ is a minimiser in \eqref{eq:H} for $\rho=\rhoa$, the first $k$ environment variables are sampled from the size-biased law $\nu'(d\zeta)=p_1(\zeta)\nu(d\zeta)/\E[p_1(\xi_1)]$ and the others from $\nu$, and $P_\xi(M_n\le0)$ is computed exactly for each sample. With $W=\prod_{i\le k}\E[p_1(\xi_1)]/p_1(\xi_i)$, the average of $W P_\xi(M_n\le0)$ is an unbiased estimator of $\E[P_\xi(M_n\le0);\,p_1(\xi_i)>0\ \forall i\le k]$. In example (A) below, $p_1>0$ in both environment states and this is $\PP(M_n\le0)$ itself; in example (B) it is the probability of a sub-event, hence a lower bound for $\PP(M_n\le0)$. In (A), the weighted contributions are highly concentrated in some runs of $1024$ samples. A bootstrap based on one such run may miss rare environments that carry substantial weight. We therefore run $8$ independent estimators of $1024$ samples each (seed $20261125+1000e+10i+j$ for example $e\in\{0,1\}$, the $i$-th value of $n$ and run $j$), report the estimate pooled over the $8192$ samples, and indicate the uncertainty by the range of the $8$ run estimates. As a check at $n=16$, where all quantities can be computed by enumeration, a single run of $1024$ samples gives for the logarithm of the target $-6.659\pm0.115$ (bootstrap standard error for this run) in (A), against the exact $\log\PP(M_{16}\le0)=-6.733$, and $-14.520\pm0.004$ in (B), against the exact logarithm $-14.518$ of the sub-event probability ($\log\PP(M_{16}\le0)=-13.055$).

\subsection{Examples}
In all examples $\mu$ is uniform on $\{-1,0,1\}$, so $L=1$, and the environment takes two values $\zeta_1,\zeta_2$ with probability $\frac12$ each:
\begin{itemize}
\item[(A)] $\zeta_1(1)=0.9$, $\zeta_1(11)=0.1$; $\zeta_2(1)=0.1$, $\zeta_2(2)=0.8$, $\zeta_2(3)=0.1$ (Schr\"oder type, constant mean $2$);
\item[(B)] $\zeta_1(1)=\zeta_1(3)=\frac12$; $\zeta_2(2)=1$ (mixed, constant mean $2$);
\item[(C)] $\zeta_1(1)=\zeta_1(2)=\frac12$; $\zeta_2(1)=0.2$, $\zeta_2(3)=0.8$ (Schr\"oder type, random mean).
\end{itemize}
For (A) and (B), \eqref{eq:constmean} holds, so Theorem~\ref{thm:annealed}(ii) applies; for (C) only Theorem~\ref{thm:annealed}(i) does. The constants are listed in Tables~\ref{tab:constants} and~\ref{tab:rates}. Quenched computations use five environment samples per example, generated with the seeds $11,\dots,15$ (\texttt{numpy.random.default_rng}); the same seed gives the same sequence of type indices in (A), (B) and (C), so the samples are coupled across examples. The single importance-sampling run used for the check at $n=16$ uses the seed $20261124$; the eight independent runs use the seeds listed in the corresponding result file.

\begin{table}[ht]
\centering
\caption{Constants of the three examples}\label{tab:constants}
\begin{tabular}{lccccc}
\toprule
 & $\ell$ & $x^*$ & $\rhoq$ & $\rhoa$ & $\beta$\\
\midrule
(A) & $0.6931$ & $0.8682$ & $1.2040$ & $0.6931$ & $0$\\
(B) & $0.6931$ & $0.8682$ & $\infty$ & $1.3863$ & $0.3466$\\
(C) & $0.6805$ & $0.8622$ & $1.1513$ & $1.0498$ & $0$\\
\bottomrule
\end{tabular}
\end{table}

\begin{table}[ht]
\centering
\caption{Quenched rates $H(x;\rhoq)$ (Theorem~\ref{thm:A}; $\infty$ means faster than exponential decay), double-exponential rates $\beta\tau^*(x)$ with $\tau^*(x)=(x^*-x)/(x^*+L)$ (Theorem~\ref{thm:B}), and $H(x;\rhoa)$, the annealed rate in (A) and (B) (Theorem~\ref{thm:annealed}(ii)) and an upper bound for $\limsup_n-\frac1n\log\PP(M_n\le xn)$ in (C) (Theorem~\ref{thm:annealed}(i); whether the limit exists is open), at $x=0$ and $x=0.3$}\label{tab:rates}
\begin{tabular}{lcccccc}
\toprule
 & \multicolumn{3}{c}{$x=0$} & \multicolumn{3}{c}{$x=0.3$}\\
\cmidrule(lr){2-4}\cmidrule(lr){5-7}
 & $H(x;\rhoq)$ & $\beta\tau^*$ & $H(x;\rhoa)$ & $H(x;\rhoq)$ & $\beta\tau^*$ & $H(x;\rhoa)$\\
\midrule
(A) & $0.8843$ & --- & $0.5593$ & $0.5787$ & --- & $0.3661$\\
(B) & $\infty$ & $0.1611$ & $0.9915$ & $\infty$ & $0.1054$ & $0.6489$\\
(C) & $0.8502$ & --- & $0.7884$ & $0.5544$ & --- & $0.5141$\\
\bottomrule
\end{tabular}
\end{table}

\subsection{Results}
\begin{figure}[ht]
\centering
\includegraphics{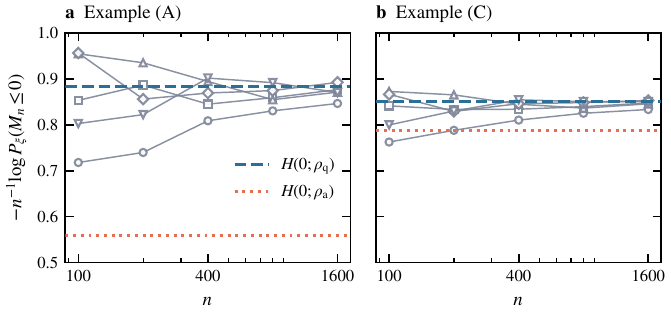}
\caption{Quenched exponents $-\frac1n\log P_\xi(M_n\le0)$ for five environment samples at $n=100,200,400,800,1600$ (logarithmic axis; lines are drawn for readability): \textbf{a} example (A), \textbf{b} example (C), $x=0$. The marker identifies the seed (circle, square, triangle, diamond, inverted triangle for seeds $11,\dots,15$); equal markers in \textbf{a} and \textbf{b} correspond to the same sequence of environment types. Dashed line: $H(0;\rhoq)$ (Theorem~\ref{thm:A}). Dotted line: $H(0;\rhoa)$, the annealed rate for (A) (Theorem~\ref{thm:annealed}(ii)) and an upper bound for $\limsup_n-\frac1n\log\PP(M_n\le0)$ in (C) (Theorem~\ref{thm:annealed}(i))}\label{fig:A}
\end{figure}

\emph{Quenched, Schr\"oder type.} Fig.~\ref{fig:A} shows $-\frac1n\log P_\xi(M_n\le0)$ for $100\le n\le1600$. At $n=1600$ the values lie in $[0.8466,0.8927]$ for (A) and in $[0.8335,0.8532]$ for (C), around the rates $0.8843$ and $0.8502$ of Theorem~\ref{thm:A}; for $x=0.3$ (not plotted) they lie in $[0.5489,0.5793]$ and $[0.5396,0.5574]$, around $0.5787$ and $0.5544$. The spread between samples is consistent with fluctuations of the environment: let
\[
\mathcal R_n(\xi;x):=\min_{1\le k\le n}\frac1n\Bigl[\sum_{i\le k}\bigl(-\log p_1(\xi_i)\bigr)+k\,\Imin\Bigl(\frac{xn-(n-k)v_k}{k}\Bigr)\Bigr],
\]
the cost of the single-particle strategy along the given environment, where $v_k$ is the speed \eqref{eq:xstar} computed with $\ell$ replaced by the empirical mean $\ell_k$ of $\log m(\xi_i)$ over $k<i\le n$ (in (A), $v_k=x^*$; for $k=n$ the second term is $n\Imin(x)$ and $v_n$ is not needed). At $n=1600$, $|{-\frac1n\log P_\xi(M_n\le xn)}-\mathcal R_n(\xi;x)|<0.0043$ for every sample and $x\in\{0,0.3\}$, compared with deviations up to $0.07$ at $n=100$, and the difference is negative in all $20$ cases. (With $v_k=x^*$ in (C), the largest deviation is $0.0049$ and two deviations are positive, both for one environment sample.) The negative sign is consistent with a polynomial prefactor larger than one, as in the constant environment, where the fitted $\hat a$ is negative.

\begin{figure}[ht]
\centering
\includegraphics{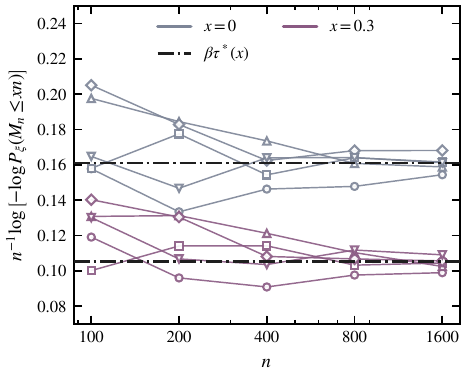}
\caption{Example (B): $\frac1n\log[-\log P_\xi(M_n\le xn)]$ for five environment samples (markers as in Fig.~\ref{fig:A}) at $n=100,\dots,1600$ (logarithmic axis), for $x=0$ (upper group of curves) and $x=0.3$ (lower group). Dash-dotted lines: $\beta\tau^*(x)$ for the two values of $x$ (Theorem~\ref{thm:B}). Markers of different seeds may coincide (e.g.\ seeds $11$ and $12$ at $n=100$, $x=0$; seeds $13$ and $15$ at $n=100$, $x=0.3$)}\label{fig:B}
\end{figure}

\emph{Quenched, mixed environment.} Fig.~\ref{fig:B} shows $\frac1n\log[-\log P_\xi(M_n\le xn)]$ in example (B). At $n=1600$ the values lie in $[0.1545,0.1682]$ for $x=0$ and in $[0.0990,0.1091]$ for $x=0.3$, compared with $\beta\tau^*=0.1611$ and $0.1054$ from Theorem~\ref{thm:B}. At $n=1600$, for each sample the value is within $0.0012$ of $\frac1n\sum_{i\le\lfloor\tau^*n\rfloor}\log b(\xi_i)$, so the deviations from $\beta\tau^*$ at this $n$ are accounted for by the fluctuations of $\sum_{i\le\tau^*n}\log b(\xi_i)$.

\begin{figure}[ht]
\centering
\includegraphics{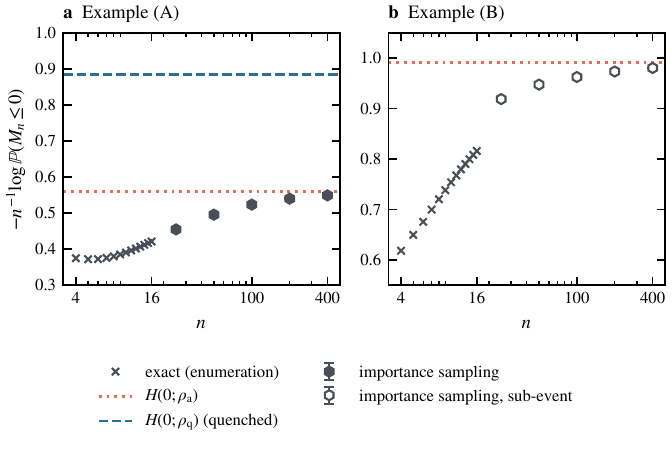}
\caption{Annealed exponents $-\frac1n\log\PP(M_n\le0)$, $x=0$: \textbf{a} example (A), \textbf{b} example (B). Crosses: exact values by enumeration of all environment sequences ($4\le n\le16$). Hexagons: importance-sampling estimates pooled over $8$ independent runs of $1024$ samples ($n=25,\dots,400$; error bars span the $8$ run estimates and are smaller than the markers); open hexagons in \textbf{b} target a sub-event and therefore give, in expectation, upper estimates of $-\frac1n\log\PP(M_n\le0)$; this explains the jump between $n=16$ and $n=25$ in \textbf{b}: at $n=16$ the sub-event gives $14.518/16=0.907$, against $0.816$ for the full probability. The $y$-ranges differ between the panels. In (B) the quenched exponent is of order $e^{0.16n}/n$ and lies far outside the plotted range. Dotted line: $H(0;\rhoa)$ (Theorem~\ref{thm:annealed}(ii)). Dashed line, \textbf{a} only: the quenched rate $H(0;\rhoq)$. Line styles for $H(0;\rhoq)$ and $H(0;\rhoa)$ are the same as in Fig.~\ref{fig:A}}\label{fig:annealed}
\end{figure}

\emph{Annealed.} Fig.~\ref{fig:annealed} shows $-\frac1n\log\PP(M_n\le0)$. The pooled importance-sampling estimates increase towards $H(0;\rhoa)$: for (A) they equal $0.454$, $0.495$, $0.523$, $0.540$, $0.549$ at $n=25,50,100,200,400$, against $H(0;\rhoa)=0.5593$; for (B) they equal $0.919$, $0.948$, $0.963$, $0.974$, $0.981$, against $0.9915$. On this scale the $8$ independent run estimates differ by at most $0.009$ in (A) (about $0.003$ for $n\ge200$) and by at most $0.001$ in (B); in (A) the effective sample sizes of single runs range from $8$ to $228$ and the largest normalised weight reaches $0.35$; in (B) they range from $26$ to $985$, with largest normalised weight $0.17$ (at $n=400$). In (B) the estimates are upper estimates of $-\frac1n\log\PP(M_n\le0)$ in expectation, since they target a sub-event; in both examples, taking the logarithm of a sample mean adds an upward bias to $-\frac1n\log$ of the estimate. In (A) the annealed values stay well below the quenched rate $H(0;\rhoq)=0.8843$, as in Corollary~\ref{cor:dichotomy}(i). In (B) the quenched probabilities of Fig.~\ref{fig:B} are of order $\exp(-e^{0.16n})$, whereas the annealed probability is of order $e^{-0.99n}$, as in Corollary~\ref{cor:dichotomy}(ii).

\FloatBarrier
\section{Discussion of the assumptions}\label{sec:discussion}

\emph{No extinction.} Assumption (A1) excludes $p_0>0$. In a constant environment, Gantert and H\"ofelsauer \cite{GantertHofelsauer2018} handle $p_0>0$ by conditioning on survival, which replaces $-\log p_1$ by $-\log\E[Z_1q^{Z_1-1}]$ with $q$ the extinction probability. In a random environment the quenched survival probability of the subtrees depends on the shifted environment, and we have not pursued this extension.

\emph{Step law.} We assume the step law does not depend on the environment and is independent of the number of children, as in \cite{ChenHe2020}. This independence is used in the many-to-one formula \eqref{eq:many} and in the first-born line. Environment-dependent step laws change $x^*$ (see \cite{MalleinMilos2019}) and, in the lower deviation strategy, the cost of moving the single particle becomes a time-inhomogeneous large deviation problem; Theorem~\ref{thm:A} is not claimed in that setting. Theorem~\ref{thm:B} uses $L<\infty$; for unbounded steps in a constant environment, Chen and He \cite[Theorems~1.3 and~1.5]{ChenHe2020} show that the decay depends on the left tail of $\mu$.

\emph{Exponential moments.} Assumption (A3) is used for the growth of $\Imin$ at $-\infty$ when $L=\infty$ (Lemmas~\ref{lem:H} and~\ref{lem:apriori}), for Cram\'er's theorem, and for $\Erw[e^{-2X}]<\infty$ in Proposition~\ref{prop:annealedupper}. We have not examined how far it can be weakened.

\emph{The annealed rate.} Without \eqref{eq:constmean}, the environment can also lower the growth rate of the population during the last $(1-t)n$ generations at an exponential cost, which slows down the maximum. By Theorem~\ref{thm:annealed}(i) and Proposition~\ref{prop:annealedupper}, the annealed probability still decays exponentially (when $\PP_{\mathrm{env}}(p_1(\xi_1)>0)>0$), but we do not know whether the exponential rate exists. We expect it to be strictly smaller than $H(x;\rhoa)$ in general, and to be given by a variational problem involving the large deviations of the environment, in analogy with the lower large deviations of branching processes in random environment \cite{BansayeBoinghoff2013}; we have not proved this. Theorem~\ref{thm:annealed}(i) remains a valid lower bound in general, so the comparison in Corollary~\ref{cor:dichotomy} does not require \eqref{eq:constmean}.

\section*{Statements and Declarations}
\paragraph*{Funding} This work received no project funding.
\paragraph*{Competing interests} The author has no competing interests to declare that are relevant to the content of this article.
\paragraph*{Author contributions} Wanting Hou is the sole author of this article.
\paragraph*{Use of generative AI} Claude (Anthropic) was used to write the numerical code and to prepare the first draft of the manuscript. The author takes full responsibility for the content of the article.
\paragraph*{Data and code availability} Online Resource~1 contains the numerical data, scripts, lattice recursion, and instructions for reproducing the tables and plots.
The numerical results were produced with Python~3.9.6, NumPy~2.0.2, SciPy~1.13.1, Matplotlib~3.9.4 and mpmath~1.3.0 by the scripts \texttt{envs.py} (definition of the examples), \texttt{run_quenched.py}, \texttt{run_annealed.py}, \texttt{run_annealed_is.py}, \texttt{run_annealed_is_reps.py}, \texttt{run_constant_check.py}, \texttt{run_validation.py} (validation and Monte Carlo check), \texttt{print_tables.py} (Tables~\ref{tab:constants} and~\ref{tab:rates}) and \texttt{make_figures.py}, which call the module \texttt{stochlib.branching.brw_lattice} implementing the recursion of Section~\ref{sec:numerics}; a README file gives the order in which the scripts are run. All computations are deterministic given the seeds stated in Section~\ref{sec:numerics} (the Monte Carlo check uses the seed $6060$, and the independent importance-sampling runs use the seeds listed in the corresponding result file), and the code reproduces all figures and tables.

\end{document}